\documentclass[12pt,oneside,english]{amsart}
\usepackage[latin1]{inputenc}
\usepackage{amssymb}

\makeatletter
\newtheorem{theorem}{Theorem}
\newtheorem{lemma}{Lemma}

\newtheorem{definition}{Definition}
\newtheorem{remark}{Remark}

\newtheorem{corollary}{Corollary}

\usepackage{babel}

\makeatother
\begin{document}
\baselineskip=17pt

\title{An Upper Estimate for the  Overpseudoprime Counting Function }

\author{Vladimir Shevelev}
\address{Department of Mathematics \\Ben-Gurion University of the
 Negev\\Beer-Sheva 84105, Israel. e-mail:shevelev@bgu.ac.il}


\begin{abstract}
 We prove that the number of overpseudoprimes to base 2 not exceeding $x$ does not exceed $x^{\frac 3 4}(1+o(1) ).$
\end{abstract}

\maketitle

\section{Introduction}

 For an odd $n >1$, consider the number $r=r(n)$ of distinct cyclotomic cosets of 2 modulo $n$  [2, pp.104-105]. E.g., $r(15)=4$
since for $n=15$ we have the following 4 cyclotomic cosets of 2: $\{1,2,4,8\}, \{3,6,12,9\},\{5,10\},\newline\{7,14,13,11\}$.
Note that, if $C_1,\ldots, C_r$ are all different cyclotomic cosets of 2$\mod n$, then

\begin{equation}\label{2}
\bigcup^r_{j=1}C_j=\{1,2,\ldots, n-1\},\qquad C_{j_1}\cap C_{j_2}=\varnothing , \;\; j_1\neq j_2.
\end{equation}

For the least common  multiple of $|C_1|, \ldots, |C_r|$ we have

\begin{equation}\label{3}
[|C_1|,\ldots,|C_r|]= h(n),
\end {equation}
where $h(n)$ is the multiplicative order of 2 modulo $n.$
(This follows easily, e.g., from Exercise 3, p. 104 in \cite{4}).

It is easy to see that for odd prime $p$ we have

\begin{equation}\label{4}
|C_1|=\ldots=|C_r|
\end {equation}

such that

\begin{equation}\label{5}
p= rh + 1.
\end {equation}

\begin{definition}
We call odd composite number $n$ overpseudoprime to base 2 $(n\in\mathbb{S}_2)$ if
\newpage
\begin{equation}\label{6}
n=r(n) h(n)+1.
\end {equation}

\end{definition}

Let $n$ be odd composite number with the prime factorization
\begin{equation}\label{7}
n=p_1^{l_1}\cdots p_k^{l_k}.
\end {equation}
In [6] we proved the following criterion.
\begin{theorem}\label{t1}
The number $n$ is overpseudoprime if and only if for all nonzero vectors $(i_1, \ldots, i_k)\leq (l_1, \ldots, l_k)$ we have
\begin{equation}\label{8}
h(n)=h(p_1^{i_1}\cdots p_k^{i_k}).
\end {equation}
\end{theorem}
\begin{corollary}\label{1} Every two overpseudoprimes $n_1$ and $n_2$ for which $h(n_1)\neq h(n_2)$ are
coprimes.
 \end{corollary}\label{1}
Notice that, every overpseudoprime is always a super-Poulet pseudoprime  and, moreover, a
strong pseudoprime to base 2 ( see Theorem 12 in [6]).
Besides,in [6] we proved the following result.
\begin{theorem}\label{2} If $p$ is prime then  $2^p-1$ is  either prime or overpseudoprime.

 \end{theorem}\label{2}
 Note that, prime divisors of overpseudoprime $n$  are primitive divisors of $2^{h(n)}-1.$ Besides, up to $2^n-1,$ every prime $p\leq n$ has already been a primitive divisor of the sequence $(2^n-1)_{n\geq1}.$ On the other hand, large prime $p>2^n-1$ evidently has $h(p)>\log_2(2p)\geq n$. Thus, in any case, all overpseudoprimes to base $2$ are in the set of products of the primitive divisors of the sequence $(2^n-1)_{n\geq1}.$ It is a simple key for finding an upper estimate for the  overpseudoprime counting function. Let $n$ be a composite number and the number $2^{n}-1$ has  at least one primitive prime divisor. Let us consider the so-called primover cofactor $([6])$ of $2^{n}-1,$ denoted $Pr(2^{n}-1)$, i. e. the products of all its primitive prime divisors. In [7] we proved that there exists a positive constant $C$ such that
\begin{equation}\label{9}
  2^{n}-1\leq (Pr(2^{n}-1))^{C\ln\ln n} .
\end{equation}
Notice also that, if $n$ is prime then, by Theorem 2, $2^{n}-1=Pr(2^{n}-1).$
\begin{definition}If $Pr(2^{n}-1)$ is not prime, then we call it full overpseudoprime to base $2$.
\end{definition}\newpage
\section{Proof of the main results}
Denote by $\omega(N(n))$ the number of prime divisors (with their multiplicities) of full overpseudoprime $N=Pr(2^{n}-1).$
\begin{lemma}\label{1}For $n>1$ we have
$$\omega(N(n))<\frac {n} {\log_2 n}.$$
 \end{lemma}
\bfseries Proof.\enskip\mdseries If $p$ is a prime divisor of $N$ then $n$ divides $p-1$ and , consequently, $p>n.$
Thus,
$$N>n^{\omega(N)}$$
and
$$\omega(N)<\frac {\log_2N} {\log_2n}<\frac {n} {\log_2n}.\blacksquare$$
Denote by $Ov^{(n)}(x)\enskip (Ov^{(\leq n)}(x))$ the number of overpseudoprimes $m\leq x$ for which $h(m)=n \enskip(h(m)\leq n).$
    \begin{lemma}\label{2}If $n\leq\log_2x,$ then
$$Ov^{(\leq n)}(x)=o(x^\varepsilon),$$
where $\varepsilon>0$ is arbitrary small for sufficiently large $x.$
\end{lemma}
\bfseries Proof.\enskip\mdseries According to Lemma 1 we, evidently, have
$$\log_2Ov^{(n)}(x)<\frac {n} {\log_2 n}\leq \frac {\log_2x} {\log_2\log_2x}.$$
Thus,
$$(Ov^{(\leq n)}(x))<x^{\frac {1} {\log_2\log_2x}}\log_2x=o(x^\varepsilon).\blacksquare $$

\begin{lemma}\label{3}If $m\leq x$ is overpseudoprime, then
$$ h(m)<x^{\frac {1} {k}}$$
and
$$k=\omega(m)\leq\frac {\log_2x} {\log_2\log_2x}.$$
\end{lemma}
\bfseries Proof.\enskip\mdseries Let $p_1\leq...\leq p_k$ be all prime divisors of overpseudoprime $m\leq x.$ Then
$$\min(p_1, ..., p_k)\leq x^{\frac {1} {k}}.$$
Thus,
$$ h(m)=h(p_1)=...=h(p_k)<x^{\frac {1} {k}}.$$
Furthermore, by Lemma 1, for $n=h(m)$ we have\newpage
$$k=\omega(m)\leq \omega(N)\leq\frac {h(m)} {\log_2h(m)}\leq\frac {x^{\frac {1} {k}}} {\log_2(x^{\frac {1} {k}})}.$$
Thus,
$$x^{\frac {1} {k}}\geq\log_2 x $$
and the lemma follows.$\blacksquare$
\begin{corollary}If $m\leq x$ is overpseudoprime, then
 $$h(m)\leq\sqrt x.$$
 \end{corollary}
Hence, denoting $Ov(x)$ the number of overpseudoprimes to base 2 not exceeding $x$, we have
$$Ov(x)=Ov^{(\leq\sqrt x)}(x).$$
\begin{lemma}\label{4} The number of overpseudoprimes $m\leq x,$ for which
$$\omega(m)=2$$
and
$$ x^{\frac {1} {4 }}\leq h(m)\leq\sqrt x,$$
is $o(x^{\frac 3 4}).$
\end{lemma}
\bfseries Proof.\enskip\mdseries
We use the following well known statement which belongs to Titchmarsh (sf [3,Theorem 5.2.1]): denote $\pi(x, k, l)$ the number of primes of the form kt+l not exceeding $x;$ if
$$ 1\leq k\leq x^{a},\enskip 0<a<1,$$
then there exists a constant $C=C(a)$ such that
 $$\pi(x, k, l)<C\frac {x} {\varphi(k)\ln x}.$$
If overpseudoprime $m=pq$ then primes $p, q$ have the form $h(m)t+1.$ Therefore, the considered number does not exceed
$$\sum_{p\leq q, pq\leq x}1=\sum_{p\leq\sqrt x}\pi(\frac {x} {p}, h(m), 1)\leq C\sum_{p\leq\sqrt x}\frac {x} {p\varphi(h(m))\ln (\frac {x} {p})}\leq$$
$$\frac {C_1x\ln\ln h(m)} {h(m)\ln \sqrt x}\sum_{p\leq\sqrt x}\frac {1} {p}\leq C_2\frac{x(\ln\ln \sqrt x)^2} {x^{\frac 1 4}\ln x}$$
and the lemma follows.$\blacksquare$\newpage
\begin{theorem}\label{3}$$Ov(x)\leq x^{\frac 3 4}(1+o(1)).$$
 \end{theorem}
\bfseries Proof.\enskip\mdseries Let $m\leq x$ be an overpseudoprime. Using the idea of C. Pomerance (private correspondence) we distinguish two cases:
a) $h(m)\leq x^{\frac {1} {2k}}$ and b) $h(m)> x^{\frac {1} {2k}}.$\newline
a) In view of Lemma 2 we could suppose that $\log_2x\leq h(m)\leq x^{\frac {1} {2k}}.$ Notice that, by Lemma 1, the number of overpseudoprimes $m\leq x,$ having $k$ prime divisors, for which $ h(m)=n$ does not exceed
$$\begin{pmatrix} \omega(N(n))\\k\end{pmatrix}\leq(\frac {h(m)} {\log_2h(m)})^{k}  $$
Summing this over $h=h(m),$ we have
$$\sum_{h=\log_2 x}^{x^\frac {1} {2k}}(\frac {h} {\log_2h})^{k}\leq\sum_{h=\log_2 x}^{x^\frac {1} {2k}}h^{k}<x^{\frac {k+1} {2k}}.$$
Further, summing this over $k\geq 2$ and using Lemma 3, we find an upper estimate of the overpseudoprimes in this case:
$$\sum_{k=2}^{\frac {\log_2 x} {\log_2\log_2 x}}x^{\frac{k+1} {2k}}\leq x^{\frac 3 4}+ x^{\frac 2 3}\frac {\log_2 x} {\log_2\log_2 x}.$$
b)In this case, using Brun-Titchmarsh inequality, for the number  of overpseudoprimes $m\leq x$ with $k\geq 3$ prime divisors we have
$$\sum_{overpseudoprimes \enskip m\leq x}\frac {1} {x}\leq \sum_{overpseudoprimes\enskip m\leq x}\frac {1} {m}\leq$$
$$ \leq (\sum_{p\leq x,  p \equiv1 \pmod h}^{k}\frac {1} {p})^{k}\leq (C_3\frac {\ln\ln x} {\varphi(h)})^{k}\leq (C_4 \frac {(\ln\ln  x)^2} {h})^{k}.$$
 Put $h_k=\max(x^{\frac {1} {2k}},\log_2 x).$ Notice that, for $k\geq 3$
$$\sum_{h_k\leq h\leq\sqrt x}\frac {1} {h^{k}}\leq \frac {C_5} {h_{k}^{k-1}}\leq \frac {C_5} {x^{\frac {k-1} {2k}}}.$$
Thus, for $k\geq3$ we have
$$\sum_{overpseudoprimes \enskip m\leq x}1\leq C_5 (C_4(\ln\ln x)^2)^{k}x^{\frac {k+1} {2k}}.$$
In view of Lemma 3,
$$(C_4(\ln\ln x)^2)^{k}=o(x^\varepsilon),$$\newpage
where $\varepsilon>0$ is arbitrary small for sufficiently large $x.$ Taking into account Lemma 4, we obtain that the number of overpseudoprimes $m\leq x$ in Case 2 for $x>x_0$ does not exceed
$$ o(x^{\frac 3 4})+ C_5x^{\varepsilon}\sum_{k\geq3}^{\frac {\log_2x} {\log_2\log_2x}}x^{\frac{k+1} {2k}}\leq o(x^{\frac 3 4})+C_5x^{\frac 2 3+\varepsilon}\frac{\log_2x} {\log_2\log_2x}=o(x^{\frac 3 4}).$$
Now, summing the numbers of overpseudoprimes $m\leq x$ in Cases 1-2, we obtain the theorem. $\blacksquare$
\begin{remark}From proof of Lemma 4, more exactly, we have
$$ Ov(x)\leq x^{\frac 3 4}(1+O(\frac {(\ln\ln x)^2} {\ln x})).$$
\end{remark}
\begin{remark} Since up to now the remainder term $O(x^{\frac 3 4})$ in the theorem on primes is unattainable, then the prime account function and the primover  account function have at the moment the same asymptotics, including the remainder term. On the other hand, C. Pomerance conjectures that really $Ov(x)=o(x^{\frac 1 2 +\varepsilon}).$ Thus, the situation, probably, will be
without changes even after proof of the Riemann hypothesis about zeros of zeta-function $(sf [3,(6.5.12)])$ .
\end{remark}
Let $Str_a(x)$ denote the number of strong pseudoprimes to base $a$ not exceeding $x.$ From Theorem 4 of paper [1] it follows that at least
\begin{equation}\label{10}
  Str_a(x)>x^{0.12-\varepsilon}.
\end{equation}
On the other hand, $Str_a(x)$ is essentially larger than $Ov_a(x).$ Indeed, for strong pseudoprime $m$ should satisfy
only conditions: $a^{m-1}\equiv1\pmod m$ and if primes $p_i|m$ then $h_a(p_i)$ contain 2 in the same powers (see [1,Proposition 1.1]).
It is interesting that (9) was  obtained in [1] for those Carmichael pseudoprimes which are strong pseudoprimes to base $a\leq e^{c_\delta(\ln\ln x)^{(1-\delta)}}$ with any fixed $\delta, 0<\delta<1,$ and  the constant $c_\delta$ depends on $\delta$ only. Recently, we have found the first Carmichael pseudoprime which is also overpseudoprime to
base 2. It is $1541955409=499*1163*2657$ such that $ h_2(499)=h_2(1163)=h_2(2657)=166.$ But it is not overpseudoprime to base 3.\newpage
\section{On overpseudoprime witness for odd composites}
For an odd composite number $n,$ let $w^{(o)}(n)$ denote the least overpseudoprime witness for $n;$ that is, the least
positive number $w^{(o)}$ for which $n$ is not an overpseudoprime to the base $w^{(o)}.$ It is very interesting to get an answer to the following Lenstra-like question: whether, for any given finite set of odd composite numbers, there exist an integer $w^{(o)}$ which serves as a witness for every number in the set (in particular, we would like  to have such a common witness for the set of odd composites up to $x.)$ Notice that, the original Lenstra's question for strong pseudoprimes was solved in [1] in negative.

   \section{Unconditional proof of infinity of overpseudoprimes to base $2$}
The following theorem belongs to C. Pomerance (private correspondence).
\begin{theorem}\label{4}
There exist infinitely many overpseudoprimes to base $2.$
\end{theorem}
\bfseries Proof.\enskip\mdseries Let $n=8k+4$. Then all primitive divisors of $2^{n}-1$ devide $2^{4k+2}+1.$ We have
the following Aurefeuillian decomposition:
$$ 2^{4k+2}+1= 4(2^{2k})^2+1=(2^{2k+1}+2^{k+1}+1)(2^{2k+1}-2^{k+1}+1)$$
and, according to [5], for every $k\geq 3$ each expression in brackets has at least one primitive divisor. Since
the difference of these expressions is a power of 2, then we have at least two different primitive divisors, for which the multiplicative order of 2 equals to $n$. Thus, product of these primitive divisors is overpseudoprime to base 2.$\blacksquare$\newline
\begin{corollary}\label{2} There exist infinitely many super-Poulet pseudoprimes to base {2}.
\end{corollary}
       So, for $n= 28, 36, 44, 52, 60, 68, 76, 84, 92, 100, 108, ...$ we have the following least overpseudoprimes to base $2$ correspondingly, with the multiplicative order of 2 which equals to $n$ (cf [8, A141232 and A122929]):
$$3277, 4033, 838861, 85489, 80581, 130561, 104653, 20647621, 280601,$$
$$818201, 68719214593, ...$$

\newpage

\bfseries Acknowledgment.\mdseries\enskip  The author is grateful to  Professor C. Pomerance  for important private correspondences.

\end{document}